\documentclass[12pt,reqno]{article}

\usepackage[usenames]{color}
\usepackage{amssymb}
\usepackage{graphicx}
\usepackage{amscd}

\usepackage[colorlinks=true,
linkcolor=webgreen,
filecolor=webbrown,
citecolor=webgreen]{hyperref}

\definecolor{webgreen}{rgb}{0,.5,0}
\definecolor{webbrown}{rgb}{.6,0,0}

\usepackage{color}
\usepackage{fullpage}
\usepackage{float}

\usepackage{graphics,amsmath,amssymb}
\usepackage{amsthm}
\usepackage{amsfonts}
\usepackage{latexsym}
\usepackage{epsf}

\newcommand{\seqnum}[1]{\href{http://oeis.org/#1}{\underline{#1}}}

\begin{document}

\begin{center}
\epsfxsize=4in
\end{center}

\theoremstyle{plain}
\newtheorem{theorem}{Theorem}
\newtheorem{corollary}[theorem]{Corollary}
\newtheorem{lemma}[theorem]{Lemma}
\newtheorem{proposition}[theorem]{Proposition}

\theoremstyle{definition}
\newtheorem{definition}[theorem]{Definition}
\newtheorem{example}[theorem]{Example}
\newtheorem{conjecture}[theorem]{Conjecture}

\theoremstyle{remark}
\newtheorem{remark}[theorem]{Remark}

\begin{center}
\vskip 1cm{\LARGE\bf 
A New Lower Bound for the  \\
\vskip .12in
Distinct Distance Constant
}
\vskip 1cm
\large
Raffaele Salvia \\
\href{mailto:raffaelesalvia@alice.it}{\tt raffaelesalvia@alice.it} \\
\end{center}

\vskip .2 in

\begin{abstract}
The reciprocal sum of Zhang sequence is not equal to the Distinct Distance Constant. This note introduces a $B_2$-sequence with larger reciprocal sum, and provides a more precise estimate of the reciprocal sums of Mian-Chowla sequence and Zhang sequence.
\end{abstract}

\vskip .2in

\section{Introduction}
A \emph{Sidon sequence}, also called a \emph{$B_2$-sequence},
is a sequence of positive integers 
$$a_1~<~a_2~<~a_3~<~\cdots$$ 
such that all the sums $a_i+a_j \ (i \leq j)$ are distinct.

The \emph{distinct distance constant} (DDC) is the supremum of the set of the reciprocal sums of Sidon sequences. Levine \cite{Levine} observed that
\begin{equation}
{\rm DDC} \leq \sum_{n=0}^\infty \frac{1}{1+\frac{n(n+1)}{2}} = \frac{2\pi}{\sqrt{7}}\tanh{ \left( \frac{\sqrt{7}}{2}\pi \right)} < 2.37366\,.
\end{equation}
Let $S_A$ be the reciprocal sum of the sequence $A$; that is,
$S_A= \sum_{i=1}^\infty \frac{1}{a_i}$. It is an open problem to find, if it exists, a Sidon sequence $U$ whose reciprocal sum is equal to the DDC \cite[p.~351]{Guy}. We only know from Taylor and Yovanof \cite{Taylor and Yovanof} that, if some Sidon sequence achieves the DDC, then it must begin with the values 1, 2, 4.

The $B_2$-sequence with the largest known reciprocal sum had been,
for a long time,
the one produced by the greedy algorithm, $G=\{ 1, 2, 4, 8, 13, 21, 31, 45, 66, 81, \dots \}$. The sequence $G$ is named the
\emph{Mian-Chowla sequence}. Lewis found that
\begin{equation}
2.158435 \leq S_G \leq 2.158677\,,
\end{equation}
where $S_G$ is the reciprocal sum of $G$ \cite[pp.~164--165]{Finch}.

In 1991, Zhang \cite{Zhang} found a Sidon sequence with a reciprocal sum greater than $S_A$. Zhang's sequence $Z$ is obtained by running the greedy algorithm for the first 14 terms, setting $z_{15}=229$, and then continuing
with the greedy algorithm. Zhang proved that

\begin{equation}
S_Z > 2.1597\,.
\end{equation}
The aim of this note is to exhibit a Sidon sequence $H$ such that $S_H > S_Z$.

\section{Computation of reciprocal sums}
\subsection{Preliminary considerations}
To estimate the reciprocal sum, we will use two basic properties of Sidon sequences: they are growing sequences, and their differences $a_i-a_j \ (i \geq j)$ are all distinct. That is,
\begin{equation}
a_i \geq a_j + (i-j), \ i \geq j
\end{equation}
\begin{equation}
a_i > \frac{i(i-1)}{2}\,.
\end{equation}
Therefore, if we know the values of $a_n$ for $1 \leq n \leq k$, we also know that
\begin{equation}
\sum_{n=1}^k \frac{1}{a_n} < S_A < \sum_{n=1}^k \frac{1}{a_n} + \sum_{n=k}^\infty \frac{1}{\max\left(a_k+n-k , \frac{n(n-1)}{2} \right)}\,.
\end{equation}
Further assumptions could be made on $a_i$, but for large values of $k$ they would not yield a significant improvement of the bounds.

\subsection{The reciprocal sum of the Mian-Chowla sequence}
Let $G$ be the $B_2$-sequence constructed by the greedy algorithm. The values of $g_n$ for $1 \leq n \leq 25000$, computed on my notebook, are listed in the accompanying file \texttt{MianChowla.txt}.
We get the following bounds for $S_G$:
\begin{equation}
\sum_{n=1}^{25000} \frac{1}{g_n} < S_G < \sum_{n=1}^{25000} \frac{1}{g_n} + \sum_{n=25001}^{510096} \frac{1}{g_{25000}+n-25000} + \sum_{n=510097}^\infty \frac{2}{n(n-1)}\,,
\end{equation}
i.e.,
\begin{equation}
2.15845268 < S_G < 2.15846062\,.
\end{equation}

\subsection{The reciprocal sum of the Zhang sequence}
The Zhang sequence $Z$ \cite{Zhang} is, at the moment, the known Sidon sequence with the largest reciprocal sum. In the accompanying file \texttt{Zhang.txt} there are the first 25000 terms of $Z$. Their values allow us to compute the following bounds for $S_Z$:
\begin{equation}
\sum_{n=1}^{25000} \frac{1}{z_n} < S_Z < \sum_{n=1}^{25000} \frac{1}{z_n} + \sum_{n=25001}^{510290} \frac{1}{z_{25000}+n-25000} + \sum_{n=510291}^\infty \frac{2}{n(n-1)}
\end{equation}
\begin{equation}
S_G < 2.16007769 < S_Z < 2.16008532\,.
\end{equation}

\subsection{The sequence $H$ and its reciprocal sum}
\begin{definition}

\begin{displaymath}
h_n = \begin{cases}
	1, & \text{if $n=1$;} \\
    229, & \text{if $n=15$;} \\
	962, & \text{if $n=27$;} \\
    \min \{ x | \forall i, j, k \leq n, a_i + a_j \neq a_k + x \} , & \text{otherwise.}
	\end{cases}
\end{displaymath}
\end{definition}
The first 26 terms of the sequence $H$ are the same of the Zhang seqence, the 27$^{th}$ term is 962, and from there on the values are provided by the greedy algorithm.

We will now prove that the sequence $H$ satisfies our aim, i.e., $S_H>S_Z$.

We can find the first values of $h_n$ ($1 \leq n \leq 25000$) in the accompanying file \texttt{H.txt}.  Hence we get
\begin{equation}
\sum_{n=1}^{25000} \frac{1}{h_n} < S_H < \sum_{n=1}^{25000} \frac{1}{h_n} + \sum_{n=25001}^{510140} \frac{1}{h_{25000}+n-25000} + \sum_{n=510141}^\infty \frac{2}{n(n-1)}
\end{equation}
\begin{equation}
2.16027651 < S_H < 2.16028417\,.
\end{equation}
In conclusion,
\begin{equation}
S_Z < 2.16027651 < S_H \leq {\rm DDC}\,,
\end{equation}
which was what we wanted to prove.

I found the sequence $H$ using an algorithm which proceeds as follows:
given a finite Sidon sequence $b_1, b_2, \ldots , b_n$, choose as
$b_{n+1}$ the value that would yield the largest reciprocal sum if the
sequence were continued with the greedy algorithm; then repeat with the
sequence $b_1, b_2, \ldots , b_n, b_{n+1}$. The program considered 20
candidate values in each step, and estimated the reciprocal sums with
the terms up to 64000. Giving as input the first 13 terms of the
Mian-Chowla sequence, I obtained the first 30 terms of $H$. I have
experimented with a number of other starting sequences, without finding
any other sequence with a reciprocal sum greater than $S_Z$.

The most likely place to find a $B_2$-sequence $X$ achieving a
reciprocal sum $S_X > S_H$, if it exists, is between the sequences
having the first 27 terms in common with $H$, or with the Zhang
sequence.

\bigskip
\hrule
\bigskip

\noindent 2010 {\it Mathematics Subject Classification}: Primary
05B10; Secondary 11Y55
\noindent \emph{Keywords: } Sidon sequences, Distinct difference constant, Mian-Chowla sequence, Zhang sequence

\bigskip
\hrule
\bigskip

\noindent (Concerned with sequence
\seqnum{A005282})


\begin{thebibliography}{1}
\bibitem{Levine}
Eugene Levine,  An extremal result for sum-free sequences,
\emph{J. Number Theor.} \textbf{12} (1980), 251–-257.
\bibitem{Guy}
Richard K. Guy, \emph{Unsolved Problems in Number Theory}, 3rd ed., 
Springer Science \& Business Media (2004), 351.
\bibitem{Taylor and Yovanof}
H. Taylor and G. S. Yovanof, $B_2$-sequences and the distinct difference constant,
\emph{Comput. Math. Appl.}  \textbf{39}(11) (2000), 37--42.
\bibitem{Finch}
Steven R. Finch, \emph{Mathematical constants},
Cambridge University Press (2003), 164--165.
\bibitem{Zhang}
Zhenxiang Zhang, A $B_2$-sequence with larger reciprocal sum,
\emph{Math. Comput.} \textbf{60}(202) (1993), 835--839.
\end{thebibliography}
\end{document}